\documentclass[12pt,leqno]{article}
\usepackage{amsmath,amssymb,amscd,latexsym}
\usepackage[all]{xy}
\newcommand{\A}{{\mathbb{A}}}
\newcommand{\C}{{\mathbb{C}}}

\newcommand{\F}{{\mathbb{F}}}

\newcommand{\Q}{{\mathbb{Q}}}

\newcommand{\R}{{\mathbb{R}}}
\newcommand{\uR}{\underline{\R}}
\newcommand{\Z}{{\mathbb{Z}}}

\newcommand{\abb}{\mathrm{ab}}
\newcommand{\ann}{\mathrm{an}}
\newcommand{\car}{\mathrm{char}}
\newcommand{\et}{\mathrm{\acute{e}t}}
\newcommand{\id}{\mathrm{id}}

\newcommand{\ord}{\mathrm{ord}}
\newcommand{\sgn}{\mathrm{sgn}\,}
\newcommand{\spec}{\mathrm{spec}\,}

\newcommand{\Fr}{\mathrm{Fr}}
\newcommand{\Gr}{\mathrm{Gr}}

\newcommand{\Imm}{\mathrm{Im}\,}

\newcommand{\map}{\mathrm{map}\,}

\newcommand{\SO}{\mathrm{SO}\,}
\newcommand{\Tr}{\mathrm{Tr}}

\newcommand{\Bh}{{\mathcal B}}
\newcommand{\Ch}{{\mathcal C}}

\newcommand{\Fh}{{\mathcal F}}

\newcommand{\Rh}{{\mathcal R}}

\newcommand{\Uh}{\mathcal{U}}
\newcommand{\oUh}{\overline{\Uh}}
\newcommand{\Vh}{{\mathcal V}}

\newcommand{\eo}{\mathfrak{o}}
\newcommand{\ep}{\mathfrak{p}}

\newcommand{\eX}{{\mathcal X}}
\newcommand{\oeX}{\overline{\eX}}
\newcommand{\oF}{\overline{F}}
\newcommand{\OF}{\overline{\F}}

\newcommand{\oX}{\overline{X}}

\newcommand{\om}{\overline{\mu}}
\newcommand{\ohne}{\smallsetminus}
\newcommand{\silo}{\stackrel{\sim}{\longrightarrow}}
\newcommand{\tei}{\, | \,}
\newcommand{\ent}{\;\widehat{=}\;}
\newcommand{\hullet}{\raisebox{0.05cm}{$\scriptscriptstyle \bullet$}}
\newcommand{\verk}{\mbox{\scriptsize $\,\circ\,$}}
\newcommand{\halb}{\frac{1}{2}}
\newcommand{\an}{\mbox{``}}
\newcommand{\ab}{\mbox{''}}
\newcommand{\dis}{\displaystyle}

\newtheorem{theorem}{Theorem}[section]

\newtheorem{cor}[theorem]{Corollary}

\newtheorem{punkt}[theorem]{$\!\!$}

\newenvironment{proof}{\bigskip{\sc Proof}}{\mbox{}\hfill$\Box$}
\begin{document}
\title{A note on arithmetic topology and dynamical systems\\[0.3cm]
{\normalsize \it Dedicated to Alexei Nikolaevich Parshin}} 
\author{Christopher Deninger}
\date{\ }
\maketitle

\section{Introduction}
\label{sec:1}

In the sixties Mazur and Manin pointed out intriguing analogies between prime ideals in number rings and knots in $3$-manifolds. Let us recall some of the relevant ideas. 

In many respects the spectrum of a finite field $\F_q$ behaves like a topological circle. For example its \'etale cohomology with $\Z_l$-coefficients for $l \nmid q$ is isomorphic to $\Z_l$ in degrees $0$ and $1$ and it vanishes in higher degrees. 

Artin and Verdier \cite{AV} have defined an \'etale topology on $\overline{\spec \Z} = \\
\spec \Z \cup \{ \infty \}$. In this topology $\overline{\spec \Z}$ has cohomological dimension three, up to $2$-torsion. 

The product formula
\[
\prod_{p \le \infty} |a|_p = 1 \quad \mbox{for} \; a \in \Q^* \; ,
\]
allows one to view $\overline{\spec \Z}$ as a {\it compact} space: Namely in the function field case the analogue of the product formula is equivalent to the formula
\[
\sum_x \ord_x f = 0
\]
on a {\it proper} curve $X_0 / \F_q$. 
Here $x$ runs over the closed points of $X_0$ and $f \in \F_q (X_0)^*$. 

By a theorem of Minkowski there are no non-trivial extensions of $\Q$ unramified at all places $p \le \infty$. Thus 
\[
\hat{\pi}_1 (\overline{\spec \Z}) = 0 \; .
\]
Hence, by analogy with the Poincar\'e conjecture, Mazur suggested to think of $\overline{\spec \Z}$ as an arithmetic analogue of the $3$-sphere $S^3$. Under this analogy the inclusion:
\[
\spec \F_p \hookrightarrow \overline{\spec \Z}
\]
corresponds to an embedded circle i.e. to a knot.

The analogue of the Alexander polynomial of a knot turns out to be the Iwasawa zeta function. One can make this precise using $p$-adic cohomology for schemes over $\spec \Z_{(p)}$. 

More generally, for a number field $k$ set
\[
\overline{\spec \eo_k} = \spec \eo_k \cup \{ \ep \tei \infty \}
\]
together with its Artin--Verdier \'etale topology. Via the inclusion
\[
\spec \eo_k / \ep \hookrightarrow \spec \eo_k
\]
we may imagine a prime ideal $\ep$ as being analogous to a knot in a compact $3$-manifold.

This nice analogy between number theory and three dimensional topology was further extended by Reznikov and Kapranov and baptized Arithmetic Topology. The reader may find dictionaries between the two fields in \cite{Re} and \cite{S}. Further analogies were contributed in \cite{Mo} and \cite{Ra} for example. In particular Ramachandran had the idea that the infinite primes of a number field should correspond to the ends of a non-compact manifold -- the analogue of $\spec \eo_k$. See the appendix to section two for his argument.

There is a completely different perspective which also leads to the idea of the finite primes being knots in some $3$-spaces \cite{D5}, \cite{D7}. One may compare the ``explicit formulas'' of analytic number theory for $k / \Q$ with certain dynamical Lefschetz trace formulas for complete flows respecting one-codimensional foliations. These two types of formulas bear a striking similarity if the dimension of the foliation is two and if the finite primes $\ep$ correspond to periodic orbits of length $\log N \ep$. The infinite primes correspond to the fixed points of the flow. 

In particular a phase space $X_{\overline{\spec \eo_k}}$ corresponding to $\overline{\spec \eo_k}$ would be three dimensional. After forgetting the parametrization the periodic orbits can thus be viewed as knots in a $3$-space.

Note that in earlier work we denoted the hypothetical space $X_{\overline{\spec \eo_k}}$ by $\an \overline{\spec \eo_k}\ab$. 

We now describe the first contribution of the paper. 
It appears that the (sheaf) cohomology of $X_{\overline{\spec \eo_k}}$ with constant coefficients should play the role of an arithmetic as opposed to geometric cohomology theory for $\overline{\spec \eo_k}$. In particular it would have cohomological dimension three. We wish to compare this as yet speculative theory with the $l$-adic cohomology of $\overline{\spec \eo_k}$. More precisely we compare Lefschetz numbers of certain endomorphisms on these cohomologies.

In order to do this we first calculate the Lefschetz number  of an automorphism $\sigma$ of $k$ on the Artin--Verdier \'etale cohomology of $\overline{\spec \eo_k}$.

Next we prove a generalization of Hopf's formula to a formula for the Lefschetz number of an endomorphism of a dynamical system on a manifold.

Assuming our formula applies to $(X_{\overline{\spec \eo_k}} , \phi^t)$ we then obtain an expression for the Lefschetz number of the automorphism on $H^{\hullet} (X_{\overline{\spec \eo_k}} , \R)$ induced by $\sigma$.

As it turns out, the two kinds of Lefschetz numbers agree in all cases, even when generalized to constructible sheaf coefficients. This result  was prompted by a question of B. Mazur on the significance of \'etale Euler characteristics in our dynamical picture \cite{D5}.

The second goal of our note is this: We show that if there is a natural dynamical system corresponding to $\spec \eo_k$ there should exist distinguished trajectories of the flow converging against the fixed points of the flow. This fits nicely with Ramachandran's idea of the infinite primes corresponding to ends of a manifold.

I am grateful to Alexei Parshin for a letter introducing me to arithmetic topology some years ago.
I would also like to thank B. Mazur for his interesting question and N. Ramachandran for allowing me to sketch his as yet unpublished ideas on ends and infinite places for the convenience of the reader.


\section{Lefschetz numbers}
\label{sec:3}

We first recall a version of \'etale cohomology with compact supports of arithmetic schemes that takes into account the fibres at infinity. Using this theory we reformulate the main result of \cite{D1} on $l$-adic Lefschetz numbers of $l$-adic sheaves on arithmetic schemes as a vanishing statement. This formulation was suggested by Faltings \cite{F} in his review of \cite{D1}. We also calculate $l$-adic Lefschetz numbers on arithmetically compactified schemes.

Let the scheme $\Uh / \Z$ be algebraic i.e. separated and of finite type and set $\Uh_{\infty} = \Uh^{\ann}_{\C} / G_{\R}$ where $\Uh_{\C} = \Uh \otimes_{\Z} \C$ and the Galois group $G_{\R}$ of $\R$ acts on $\Uh^{\ann}_{\C}$ by complex conjugation. We give $\Uh_{\infty}$ the quotient topology of $\Uh^{\ann}_{\C}$.

Artin and Verdier \cite{AV} define the \'etale topology on $\oUh = \Uh \amalg \Uh_{\infty}$ as follows. The category of ``open sets'' has objects the pairs \\
$(f : \Uh' \to \Uh, D')$ where $f$ is an \'etale morphism and $D' \subset \Uh'_{\infty}$ is open. The map $f_{\infty} : D' \to \Uh_{\infty}$ induced by $f$ is supposed to be ``unramified'' in the sense that $f_{\infty} (d') \in \Uh (\R)$ if and only if $d' \in \Uh' (\R)$. Note that $\Uh (\R)$ is a closed subset of $\Uh_{\infty}$.

A morphism
\[
(f : \Uh' \to \Uh , D') \longrightarrow (g : \Uh'' \to \Uh , D'')
\]
is a map $\Uh' \to \Uh''$ commuting with the structure maps and such that the induced map $\Uh'_{\infty} \to \Uh''_{\infty}$ carries $D'$ into $D''$. 

Coverings are the obvious ones.

Pullback defines morphisms of sites:
\[
\Uh_{\et} \overset{j}{\longrightarrow} \oUh_{\et} \overset{i}{\longleftarrow} \Uh_{\infty} \; .
\]
Let $\sim$ denote the corresponding categories of abelian sheaves. One proves that $\oUh^{\sim}_{\et}$ is the mapping cone of the left exact functor
\[
i^* j_* : \tilde{\Uh}_{\et} \longrightarrow \tilde{\Uh}_{\infty} \; .
\]
In particular we have maps $i^!$ and $j_!$ at our disposal. Let us describe the functor $i^* j_*$ explicitely. Let $\alpha : \Uh^{\ann}_{\C} \to \Uh_{\et}$ be the canonical map of sites. Note that $\alpha^* F$ is a $G_{\R}$-sheaf on $\Uh^{\ann}_{\C}$ for every sheaf $F$ on $\Uh_{\et}$. If $\pi : \Uh^{\ann}_{\C} \to \Uh_{\infty}$ denotes the natural projection, define the left exact functor
\[
\pi^{G_{\R}}_* : (\mbox{abelian}\;G_{\R}\mbox{-sheaves on} \; \Uh^{\ann}_{\C}) \longrightarrow \tilde{\Uh}_{\infty}
\]
by
\[
\pi^{G_{\R}}_* (G) (V) = G (\pi^{-1} (V))^{G_{\R}} \; .
\]
Then one can check that
\[
i^* j_* = \pi^{G_{\R}}_* \verk \alpha^* \; .
\]
In particular we see that
\[
i^* R^n j_* = R^n \pi^{G_{\R}}_* \verk \alpha^* \; .
\]
It follows that for $n \ge 1$ the sheaf $R^n j_* F = i_* i^* R^n j_* F$ is $2$-torsion with support on $\Uh (\R) \subset \Uh_{\infty}$ as stated in \cite{AV}. In particular
\begin{equation}
  \label{eq:11}
  H^n_{\et} (\oUh , j_* F) \longrightarrow H^n_{\et} (\oUh , R j_* F) = H^n_{\et} (\Uh , F)
\end{equation}
is an isomorphism up to $2$-torsion for $n \ge 1$ (and an isomorphism for $n = 0$).

Let us now define cohomology with compact supports. Assume first that $\eX / \Z$ is a proper scheme. We define:
\[
H^n_c (\eX , F) := H^n (\oeX , j_! F) \; .
\]
The distinguished triangle
\[
j_! j^* G^{\hullet} \longrightarrow G^{\hullet} \longrightarrow i_* i^* G^{\hullet} \longrightarrow \ldots
\]
for complexes of sheaves $G^{\hullet}$ on $\oeX$ applied to $G^{\hullet} = Rj_* F$ gives the triangle
\[
j_! F \longrightarrow Rj_* F \longrightarrow i_* i^* Rj_* F \longrightarrow \ldots
\]
From this one gets an exact sequence
\[
\longrightarrow H^n_c (\eX , F) \longrightarrow H^n (\eX , F) \longrightarrow H^n (\eX_{\infty} , i^* Rj_* F) \longrightarrow \ldots
\]
We have $i^* Rj_* = R\pi^{G_{\R}}_* \verk \alpha^*$ and
\begin{eqnarray}
H^n (\eX_{\infty} , i^* Rj_* F) & = & 
  H^n (\eX_{\infty} , R \pi^{G_{\R}}_* (\alpha^* F)) \nonumber \\
& = & H^n (G_{\R} , R\Gamma (\eX^{\ann}_{\C} , \alpha^* F)) \nonumber \\
& = & H^n (\eX^{\ann}_{\R} , \alpha^* F) \; . \label{eq:12}
\end{eqnarray}
Here, for any $G_{\R}$-sheaf $G$ on a complex analytic space $Y$ with a real structure we set:
\[
H^n (Y_{\R} , G) = H^n (G_{\R} , R \Gamma (Y , G)) \; .
\]
In conclusion we get the long exact sequence:
\[
  \longrightarrow H^n_c (\eX , F) \longrightarrow H^n (\eX, F) \longrightarrow H^n (\eX^{\ann}_{\R} , \alpha^* F) \longrightarrow \ldots
\]
We now consider the case where $\Uh / \Z$ is algebraic but not necessarily proper. Choose a Nagata compactification of $\Uh$ i.e. a proper $\Z$-scheme $\eX$ with an open immersion $\kappa : \Uh \to \eX$. One sets:
\[
H^n_c (\Uh , F) = H^n_c (\eX , \kappa_! F) \; .
\]
This group and the group
\[
H^n_{fc} (\Uh , F) = H^n (\eX , \kappa_! F)
\]
can be shown not to depend on the compactification $\eX$. Beware that in \cite{D1} the cohomology $H^n_{fc}$ was denoted $H^n_c$.

Since $\alpha^* \kappa_! = \kappa_! \alpha^*$ we obtain a long exact sequence:
\begin{equation}
  \label{eq:13}
  \longrightarrow H^n_c (\Uh , F) \longrightarrow H^n_{fc} (\Uh , F) \longrightarrow H^n_c (\Uh^{\ann}_{\R} , \alpha^* F) \longrightarrow \ldots
\end{equation}
The definition of cohomology can be extended to $\Q_l$-sheaves and the exact sequence (3) continues to hold. See \cite{D1} \S\,3 in this regard. The maps in (1) induce isomorphisms
\begin{equation}
  \label{eq:14}
  H^n (\oUh , j_* F) \silo H^n (\Uh , F) \quad \mbox{for $\Q_l$-sheaves $F$ and} \; n \ge 0 \; .
\end{equation}
An endomorphism $(\sigma , e)$ of the pair $(\Uh , F)$ is a morphism $\sigma : \Uh \to \Uh$ together with a homomorphism $e : \sigma^* F \to F$ of $\Q_l$-sheaves. It induces pullback endomorphisms on $H^n_c (\Uh , F)$ and $H^n (\Uh , F)$ and on $H^n_c (\Uh^{\ann}_{\R} , \alpha^* F)$. For an endomorphism $\varphi$ of a finite dimensional graded vector space $H^{\hullet}$ we will abbreviate the Lefschetz number as follows:
\[
\Tr (\varphi \tei H^{\hullet}) := \sum_i (-1)^i \Tr (\varphi \tei H^i) \; .
\]
The following theorem is a reformulation of results in \cite{D1}. They are based on algebraic number theory and in particular on class field theory.

\begin{theorem}
  \label{t31}
Let $(\sigma , e)$ be an endomorphism of $(\Uh , F)$ as above, $l \neq 2$. The cohomologies $H^{\hullet}_c (\Uh , F)$ and $H^{\hullet} (\Uh , F)$ are finite dimensional and we have:
\begin{equation}
  \label{eq:15}
  \Tr ((\sigma , e)^* \tei H^{\hullet}_c (\Uh , F)) = 0 \; .
\end{equation}
If $\Uh = \eX$ is proper over $\spec \Z$ then:
\begin{equation}
  \label{eq:16}
  \Tr ((\sigma , e)^* \tei H^{\hullet} (\oeX , j_* F)) = \Tr ((\sigma , e)^* \tei H^{\hullet} (\eX^{\ann}_{\R} , \alpha^* F)) \; .
\end{equation}
If in addition $\eX$ is generically smooth and the fixed points $x$ of $\sigma$ on $\eX_{\infty}$ are non-degenerate in the sense that $\det (1 - T_x\sigma \tei T_x (\eX \otimes \R))$ is non-zero then we have:
\begin{equation}
  \label{eq:17}
  \Tr ((\sigma , e)^* \tei H^{\hullet} (\oeX , j_* F)) = \sum_{x \in \eX_{\infty} \atop \sigma x = x} \Tr (e_x \tei (\alpha^* F)_x) \varepsilon_x (\sigma) \; .
\end{equation}
Here
\[
\varepsilon_x (\sigma) = \sgn \det (1 - T_x \sigma \tei T_x (\eX \otimes \R)) \; .
\]
\end{theorem}

\begin{proof}
  According to \cite{D1} (3.10) we have for $l \neq 2$:
  \begin{equation}
    \label{eq:18}
    \Tr ((\sigma , e)^* \tei H^{\hullet}_{fc} (\Uh , F)) = \Tr ((\sigma , e)^* \tei H^{\hullet}_c (\Uh \otimes \R , F)) \; .
  \end{equation}
Now for a constructible $\Z / l^n$-sheaf $F_n$ on $\Uh_{\C}$ we have
\[
R \Gamma_c (\Uh^{\ann}_{\C} , F_n) = R\Gamma_c (\Uh_{\C} , F_n)
\]
by the comparison theorem between ordinary and \'etale cohomology. If $F_n$ carries a $G_{\R}$-action relative to the $G_{\R}$-action on $\Uh_{\C}$ we get:
\[
H^{\nu}_c (\Uh^{\ann}_{\R} , F_n) = H^{\nu}  (G_{\R} , R \Gamma_c (\Uh^{\ann}_{\C} , F_n)) = H^{\nu} (G_{\R} , R\Gamma_c (\Uh_{\C} , F_n)) = H^{\nu}_c (\Uh \otimes \R , F_n) \; .
\]
Passage to the inverse limit shows that:
\begin{equation}
  \label{eq:19}
  H^{\hullet}_c (\Uh \otimes \R , F) = H^{\hullet}_c (\Uh^{\ann}_{\R} , \alpha^* F) \; .
\end{equation}
This also holds for $l = 2$.

Together with (8) and the exact sequence (3) we get the first formula (5).

The second formula follows from (4), (8) and (9).

To prove the third, note that since $F$ has no $2$-torsion:
\begin{eqnarray}
  \label{eq:20}
\Tr ((\sigma , e)^* \tei H^{\hullet} (\eX^{\ann}_{\R} , \alpha^* F)) = \Tr ((\sigma , e)^* \tei H^{\hullet} (\eX^{\ann} , \alpha^* F)^{G_{\R}}) \\
\hspace*{1cm} = |G_{\R}|^{-1} \sum_{\tau \in G_{\R}} \Tr ((\sigma \tau , \tau^* e)^* \tei H^{\hullet} (\eX^{\ann} , \alpha^* F)) \; . \nonumber
\end{eqnarray}
Here we have used the following elementary formula: Let $V$ be a finite dimensional $K$-vector space with an endomorphism $\varphi$. Assume that a finite group $G$ of order prime to $\car (K)$ operates on $V$ such that $\varphi \verk g = g \verk \varphi$ for all $g \in G$. Then we have:
\[
\Tr (\varphi \tei V^G) = |G|^{-1} \sum_{g \in G} \Tr (\varphi \verk g \tei V) \; .
\]
The proof uses the projection operator $|G|^{-1} \sum_{g \in G} g$ from $V$ to $V^G$. \\
Also note that the map induced by $\sigma$ on $\eX^{\ann} = \eX (\C)^{\ann}$ commutes with the $G_{\R}$-action on $\eX^{\ann}$ by complex conjugation and that $\alpha \verk \tau = \alpha$. 

Applying the Lefschetz fixed point formula to (10) we find
\begin{eqnarray*}
  \lefteqn{\Tr ((\sigma , e)^* \tei H^{\hullet} (\eX^{\ann}_{\R} , \alpha^* F))} \\
 & = & |G_{\R}|^{-1} \sum_{\tau \in G_{\R}} \sum_{x' \in \eX^{\ann} \atop \sigma \tau (x') = x'} \Tr ((\tau^* e)_{x'} \tei (\alpha^* F)_{x'}) \sgn \det (1 - T_{x'} (\sigma \tau) \tei T_{x'} (\eX_{\C})) \\
& = & \sum_{x \in \eX_{\infty}\atop \sigma x = x} \Tr (e_x \tei (\alpha^* F)_x) \varepsilon_x (\sigma) \; .
\end{eqnarray*}
\end{proof}

The following assertion is an immediate consequence of the theorem.

\begin{cor}
  \label{t32}
I) Let $\sigma$ be an endomorphism of an algebraic scheme $\Uh / \Z$. Then we have for any $l \neq 2$:
\begin{equation}
  \label{eq:21}
  \Tr (\sigma^* \tei H^{\hullet}_c (\Uh , \Q_l)) = 0 \; .
\end{equation}
If $\Uh = \eX$ is proper over $\spec \Z$ then
\begin{equation}
  \label{eq:22}
  \chi (H^{\hullet} (\oeX , j_* \Q_l)) = \chi (H^{\hullet} (\eX^{\ann}_{\R} , \Q_l)) \; .
\end{equation}
If in addition $\eX$ is generically smooth and the fixed points of $\sigma$ on $\eX_{\infty}$ are non-degenerate we have:
\begin{equation}
  \label{eq:23}
  \Tr (\sigma^* \tei H^{\hullet} (\oeX , j_* \Q_l)) = \sum_{x \in \eX_{\infty} \atop \sigma x = x} \varepsilon_x (\sigma) \; .
\end{equation}
II) Let $M$ be a motive over $k$ and view its $l$-adic realization $M_l$ as a $\Q_l$-sheaf on $\spec k$. Define the $\Q_l$-sheaf $F_l (M)$ on $\spec \eo_k$ by $F_l (M) = j_{0*} M_l$ where $j_0 : \spec k \hookrightarrow \spec \eo_k$ is the inclusion. Let $\sigma$ be an automorphism of $\spec \eo_k$. It is induced by an automorphism $\sigma_k$ of $k$. Let $e : M^{\sigma_k} \to M$ be a morphism of motives over $k$ and let $e : \sigma^* F_l (M) \to F_l (M)$ be the induced map of $\Q_l$-sheaves on $\spec \eo_k$. Then we have:
\begin{equation}
  \label{eq:24}
  \Tr ((\sigma , e)^* \tei H^{\hullet} (\overline{\spec \eo_k} , j_* F_l (M))) = \sum_{\ep \tei \infty \atop \sigma \ep = \ep} \Tr (e_{\ep} \tei M_{\ep}) \; .
\end{equation}
Here for $\ep$ with $\sigma \ep = \ep$ we let $e_{\ep}$ be the map induced by $e$ on the real Hodge realization $M_{\ep}$ of $M$.
\end{cor}

\begin{punkt}
\rm We now explain analogies and an interesting difference with the case of varieties over finite fields. For a variety $X / \F_q$ with an endomorphism $\sigma$ of $X$ over $\F_q$ and a constructible $\Q_l$-sheaf $F$ on $X$ with an endomorphism $e : \sigma^* F \to F$ we have
\begin{equation}
  \label{eq:25}
  \Tr ((\sigma , e)^* \tei H^{\hullet}_c (X,F)) = 0
\end{equation}
and
\begin{equation}
  \label{eq:26}
  \Tr ((\sigma , e)^* \tei H^{\hullet} (X,F)) = 0 \; .
\end{equation}
Contrary to the number field case these assertions are easy to prove. E.g. for the first one, set $\oX = X \otimes \overline{\F}_q , \oF = F \, |_{\oX}$. The Hochschild--Serre spectral sequence degenerates into short exact sequences:
\[
0 \longrightarrow H^1 (\F_q , H^{n-1}_c (\oX , \oF)) \longrightarrow H^n_c (X , F) \longrightarrow H^0 (\F_q , H^n_c (\oX , \oF)) \longrightarrow 0 \; .
\]
Moreover for every $G_{\F_q}$-module $M$ there is an exact sequence:
\[
0 \longrightarrow H^0 (\F_q , M) \longrightarrow M \xrightarrow{1-\varphi} M \longrightarrow H^1 (\F_q , M) \longrightarrow 0 
\]
where $\varphi$ is a generator of $G_{\F_q}$. This implies (15) and (16) follows similarly since all groups are known to be finite dimensional.
\end{punkt}

Now in the number field case, according to Theorem 2.1 (5) the analogue of (15) is valid. In particular, for any finite set $S$ of prime ideals in $\eo_k$ we have
\begin{equation}
  \label{eq:27}
  \chi (H^{\hullet}_c (\spec \eo_{k,S} , \Q_l)) = 0 \; .
\end{equation}
The analogue of (16) however is not valid in general: From Theorem 2.1 (6) it follows for example that
\begin{equation}
  \label{eq:28}
  \chi (H^{\hullet} (\overline{\spec \eo_k} , j_* \Q_l) = \; \mbox{number of infinite places of} \; k \; .
\end{equation}
The fact, that in (17) the Euler characteristic is unchanged if we increase $S$ i.e. take out more finite places $\ep$ follows directly from the fact that
\[
\chi (H^{\hullet} (\spec (\eo_k / \ep) , \Q_l)) = 0 \; .
\]
The topological intuition behind this equation is that finite primes are like circles (whose Euler characteristic also vanishes). The difference between (17) and (18) comes from the different nature of the infinite places. Cohomologically the complex places behave like points and the real places behave like the ``quotient of a point by $G_{\R}$''. This fits well with our dynamical considerations in the next section. For totally different reasons they suggest that a complex resp. real place gives rise to an embedded real line 
\[
\R \hookrightarrow X_{\overline{\spec \eo_k}}
\]
resp. to an embedded orbifold 
\[
\R / \mu_2 = \R^{\ge 0 } \hookrightarrow X_{\overline{\spec \eo_k}} \; .
\]
Contrast this with the idea that finite places should give rise to embeddings of circles 
\[
S^1 \hookrightarrow X_{\overline{\spec \eo_k}}\; .
\]

\begin{punkt}
  \rm In this section we prove a formula for Lefschetz numbers of finite order automorphisms of dynamical systems. We also allow certain constructible sheaves as coefficients.

Let us consider a complete flow $\phi^t$ on a compact manifold $X$. Assume that for every fixed point $x$ there is some $\delta_x > 0$ such that
\[
\det (1 - T_x \phi^t \tei T_x X) \neq 0 \quad \mbox{for} \quad 0 < t < \delta_x \; .
\]
This condition is weaker than non-degeneracy. Still it implies that the fixed points are isolated and hence finite in number.

We also assume that the lengths of the closed orbits are bounded below by some $\varepsilon > 0$.

Let $\sigma$ be an automorphism of finite order of $X$ which commutes with all $\phi^t$ i.e. an automorphism of $(X , \phi)$. Consider a constructible sheaf of $\Q$\,- or $\R$-vector spaces on $X$ together with an endomorphism $e : \sigma^{-1} F \to F$. We also assume that there is an action $\psi^t$ over $\phi^t$ i.e. isomorphisms $\psi^t : (\phi^t)^{-1} F \to F$ for all $t$ satisfying the relations
\[
\psi^0 = \id \quad \mbox{and} \quad \psi^{t_1 + t_2} = \psi^{t_2} \verk (\phi^{t_2})^{-1} (\psi^{t_1}) \quad \mbox{for all} \; t_1 , t_2 \in \R \; .
\]
For constant $F$ there is a canonical action $\psi$.\\
Note that $\sigma$ and $e$ determine an endomorphism $(\sigma , e)^*$ of $H^{\hullet} (X, F)$.\\
Let us say that a point $x \in X$ is $\phi$-fix if $\phi^t (x) = x$ for all $t \in \R$. 

Then we have the following formula:
\end{punkt}

\begin{theorem}
  $\dis \Tr ((\sigma , e)^* \tei H^{\hullet} (X, F)) = \sum_{x\in X , \phi- \mathrm{fix} \atop \sigma x = x} \Tr (e_x \tei F_x) \varepsilon_x (\sigma)$\\
where
\[
\varepsilon_x (\sigma) = \lim_{t \to 0 \atop t > 0} \sgn \det (1 - T_x (\phi^t \sigma ) \tei T_x X)\; .
\]
\end{theorem}

\begin{proof}
  Assume $N \ge 1$ is such that $\sigma^N = 1$ and fix some $s$ with \\
$0 < s < N^{-1} \min_{x\in X , \phi - \mathrm{fix}} (\varepsilon , \delta_x)$. If $x$ is any point of $X$ with $(\phi^s \sigma) (x) = x$, then $(\phi^{Ns} \sigma^N) (x) = x$ and hence $\phi^{Ns} (x) = x$. If $x$ lay on a periodic orbit $\gamma$ then $l (\gamma) \le Ns$ and hence $\varepsilon \le Ns$ contrary to the choice of $s$. Thus $x$ is a fixed point of $\phi$. Because of $\phi^s (x) = x$ we have $\sigma x = x$ as well. If $1$ is an eigenvalue of $T_x (\phi^t \sigma)$ then $1$ is an eigenvalue of $(T_x (\phi^t \sigma))^N = T_x \phi^{Nt}$ as well. Thus $Nt \ge \delta_x$. By assumption on $s$ we have $Ns < \delta_x$ and therefore
\[
\det (1 - T_x (\phi^s \sigma) \tei T_x X) \neq 0 \; .
\]
In conclusion: the fixed points of the automorphism $\phi^s \sigma$ coincide with those fixed points of the flow $\phi$ which are also kept fixed by $\sigma$. They are all non-degenerate.

Consider the morphisms for $t \in \R$
\[
e_t = e \verk \sigma^{-1} (\psi^t) : (\phi^t \sigma)^{-1} F \longrightarrow F \; .
\]
The Lefschetz fixed point formula for the endomorphism $(\phi^s \sigma , e_s)$ of $(X , F)$ for a fixed $s$ as above now gives the formula:
\[
\Tr ((\phi^s \sigma , e_s)^* \tei H^{\hullet} (X ,F)) = \sum_{x \in X ,  \phi-\mathrm{fix} \atop \sigma x = x } \Tr ((e_s)_x \tei F_x) \sgn \det (1 - T_x (\phi^s \sigma) \tei T_x X) \; .
\]
The left hand side is defined for all $s$ in $\R$ and by homotopy invariance of cohomology it is independent of $s$. Passing to the limit $s \to 0$ for positive $s$ in the formula thus gives the assertion.
\end{proof}

\begin{cor}
  Let $(X , \phi^t , \sigma)$ be as in theorem 2.5 and let $U \subset X$ be open, $\phi$- and $\sigma$-invariant. Assume that $X \ohne U$ is a compact submanifold of $X$. Consider a constructible sheaf $F$ of $\Q$- or $\R$-vector spaces on $U$ with an endomorphism $e : \sigma^{-1} F \to F$ and an action $\psi^t$ over $\phi^t \, |_U$. Then we have:
\[
\Tr ((\sigma , e)^* \tei H^{\hullet}_c (U , F)) = \sum_{x \in U , \phi-\mathrm{fix} \atop \sigma x = x} \Tr (e_x \tei F_x) \varepsilon_x (\sigma) \; .
\]
In particular
\[
\Tr ((\sigma , e)^* \tei H^{\hullet}_c (U ,F)) = 0
\]
if $\phi$ has no fixed points on $U$.
\end{cor}

\begin{proof}
  Apply 2.5 to $(X , j_! F)$, where $j : U \hookrightarrow X$ is the inclusion.
\end{proof}

\begin{punkt}
\rm  We now compare arithmetic and dynamic Lefschetz numbers. This section is of a heuristic nature: we assume that a functor $X_{\Uh} \mapsto (X_{\Uh} , \phi^t)$ from flat algebraic schemes over $\spec \Z$ to dynamical systems exists, with properties as described in \cite{D5} end of \S\,4. In that reference the space $X_{\Uh}$ was denoted by $\an \Uh \ab$.

We consider various cases.
\bigskip

{\sc 2.7.1} Let $\sigma$ be a finite order automorphism of the flat algebraic scheme $\Uh / \Z$. The induced automorphism $\sigma$ of $X_{\Uh}$ is an automorphism of $(X_{\Uh} , \phi^t)$ and hence commutes with each $\phi^t$. The phase space $X_{\Uh}$ cannot be a manifold \cite{D7} \S\,5. Still let us assume that the assertion of corollary 2.6 applies to
\[
(U , \phi^t , \sigma , F , e , \psi^t) = (X_{\Uh}, \phi^t , \sigma , \Q , \id , \id) \; .
\]
Then we find
\[
\Tr (\sigma^* \tei H^{\hullet}_c (X_{\Uh} , \Q)) = 0
\]
in accordance with formula (11) since $\phi^t$ should have no fixed points on $X_{\Uh}$.
\bigskip

{\sc 2.7.2} Now consider the case where in addition $\Uh = \eX$ is proper over $\spec \Z$ and generically smooth. Also the fixed points of $\sigma$ on $\eX_{\infty}$ should be non-degenerate. Under assumptions as in {\bf 1}, from corollary 2.6 applied to
\[
(U , \phi^t , \sigma, F, e , \psi^t) = (X_{\oeX} , \phi^t , \sigma , j_* \Q , \id , \id)
\]
we would get a formula corresponding to (13):
\begin{equation}
  \label{eq:29}
  \Tr (\sigma^* \tei H^{\hullet} (X_{\oeX} , j_* \Q)) = \sum_{x \in \eX_{\infty} \atop \sigma x = x} \varepsilon_x (\sigma) \; .
\end{equation}
Here we have used that the fixed point set of $\phi^t$ on $X_{\oeX}$ should be $\eX (\C) / G_{\R} = \eX_{\infty}$: In the case of $\eX = \spec \eo_k$, this is the set of archimedean valuations of $k$. For general $\eX$ we only have the following argument: The set of closed points of $\eX$ over $p$ can be identified with the set $\eX (\OF_p) / \langle \Fr_p \rangle$ of Frobenius orbits on $\eX (\OF_p)$. Thus the set of closed orbits of $(X_{\eX} , \phi^t)$ would be in bijection with the union of all $\eX (\OF_p) / \langle \Fr_p \rangle$. Correspondingly it looks natural to assume that the set of fixed points of $\phi$ would be in bijection with $\eX (\C) / \langle F_{\infty} \rangle = \eX_{\infty}$ where the infinite Frobenius $F_{\infty} : \eX (\C) \to \eX (\C)$ acts by complex conjugation. 

Actually the comparison of formulas (13) and (19) lends further credibility to this idea. 

The different definitions of $\varepsilon_x (\sigma)$ in (13) and (19) should agree. We will discuss this aspect in the next case.
\bigskip

{\sc 2.7.3} Let $M$ be a motive over the number field $k$. By the construction of \cite{D5} \S\,5 it should give rise to a constructible sheaf $F (M)$ of $\R$-vector spaces with an action $\psi^t$ over $(X_{\spec \eo_k} , \phi^t)$. Let $(\sigma , e)$ be an endomorphism of $(\spec \eo_k , M)$ as in corollary 2.2 II. The induced endomorphism of $((X_{\spec \eo_k} , \phi^t) , F (M))$ will also be denoted $(\sigma , e)$. Note that $\sigma$ has finite order. Let $j : X_{\spec \eo_k} \hookrightarrow X_{\overline{\spec \eo_k}}$ denote the inclusion. Then by theorem 2.5 we expect the formula:
\begin{equation}
  \label{eq:30}
  \Tr ((\sigma , e)^* \tei H^{\hullet} (X_{\overline{\spec \eo_k}} , j_* F(M)) =\!\!\sum_{\ep \tei \infty \atop \sigma (x_{\ep}) = x_{\ep}}\!\!\Tr (e_{x_{\ep}} \tei (j_* F (M))_{x_{\ep}}) \varepsilon_{x_{\ep}} (\sigma) \; .
\end{equation}
Here $x_{\ep}$ would be the fixed point of $\phi$ corresponding to the infinite prime $\ep$.

Comparing this with formula (13) in \cite{D6}, suggests that we have an isomorphism:
\[
(j_* F (M))_{x_{\ep}} = \Gr^{\hullet}_{\Vh} M_{\ep} \; .
\]
Here $\Vh$ is a certain filtration on the real Hodge realization of $M \otimes k_{\ep}$.
Provided that $\varepsilon_{x_{\ep}} (\sigma) = 1$ for all $\ep$, formula (20) thus becomes
\begin{eqnarray*}
  \Tr ((\sigma , e)^* \tei H^{\hullet} (X_{\overline{\spec \eo_k}} , j_* F (M))) & = & \sum_{\ep \tei \infty \atop \sigma \ep = \ep} \Tr (e_{\ep} \tei \Gr^{\hullet}_{\Vh} M_{\ep}) \\
& = & \sum_{\ep \tei \infty \atop \sigma \ep = \ep} \Tr (e_{\ep} \tei M_{\ep})
\end{eqnarray*}
in accordance with (14). 

We finish this subsection with an argument for $\varepsilon_{x_{\ep}} (\sigma) = 1$. Set $X = X_{\overline{\spec \eo_k}}$ with its $1$-codimensional foliation $\Bh$ for brevity and set $x = x_{\ep}$. On $T_x X / T_x \Bh$ the map $T_x \phi^t$ has eigenvalue $e^{\kappa_{\ep} t}$ where $\kappa_{\ep} = -1$ if $\ep$ is complex and $\kappa_{\ep} = -2$ if $\ep$ is real. According to \cite{D7} 5.7 5) the map $T_x \phi^t$ has the form $e^{t/2} O_t$ on $T_x \Bh$ for some $O_t \in \SO (T_x \Bh)$. Hence the eigenvalues of $T_x \phi^t$ on $T_x \Bh$ are complex conjugate. By functoriality, $\sigma$ respects $\Bh$, hence it fixes the leaf through $x$ if $\sigma (x) = x$. In this case $T_x \sigma$ respects $T_x \Bh$ and $T_x X / T_x \Bh$ therefore. On the $1$-dimensional real vector space $T_x X / T_x \Bh$ the map $T_x \sigma$ has eigenvalue $+1$ or $-1$ since $T_x \sigma$ has finite order. Hence for any $t> 0$ we find:
\[
\det (1 - T_x (\phi^t \sigma) \tei T_x X / T_x \Bh) = 1 - (\pm 1) e^{\kappa_{\ep} t} > 0 \; .
\]
On $T_x \Bh$, either $T_x \sigma$ has complex conjugate eigenvalues in which case $T_x (\phi^t \sigma)$ has complex conjugate eigenvalues so that 
\[
\det (1 - T_x (\phi^t \sigma) \tei T_x X / T_x \Bh) > 0 \; .
\]
Or $T_x \sigma$ has two different real eigenvalues which then must be $1$ and $-1$. But this case cannot occur since for small $t> 0$ the eigenvalues $\mu_t , \overline{\mu}_t$ of $T_x \phi^t$ on $T_x \Bh$ are not real. The eigenvalues of $T_x (\phi^t \sigma)$ on $T_x \Bh$ would be $\mu_t , - \overline{\mu}_t$ (or $- \mu_t, \overline{\mu}_t$) hence non-real hence complex conjugate i.e. $\om_t = - \om_t$, contradiction. 

In conclusion we get $\varepsilon_{x_{\ep}} (\sigma) = 1$ for all $\ep \tei \infty$. 
\end{punkt}

\bigskip
\noindent {\large \sc Appendix}

\bigskip

  Here we discuss ends of arithmetic varieties. The basic insight here is due to N. Ramachandran \cite{Ra2}. Let $X$ be a connected, locally connected, locally compact Hausdorff space. The space of ends of $X$ is the pro-discrete topological space
\[
E [X] = \lim_{\xleftarrow[K]{}} \pi_0 (X \ohne K) \; ,
\]
where $K$ runs over the compact subsets of $X$.

For a commutative ring $A$ let $C^{\infty} (E [X] , A)$ denote the $A$-module of locally constant $A$-valued functions on $E [X]$. Thus
\begin{eqnarray*}
  C^{\infty} (E [X] , A) & = & \lim_{\xrightarrow[K]{}} \map (\pi_0 (X \ohne K) , A) \\
& = & \lim_{\xrightarrow[K]{}} H^0 (X \ohne K , A) \; .
\end{eqnarray*}
Here as usual we use sheaf cohomology. Passing to the limit over $K$ in the exact sequence:
\[
0 \to H^0_K (X ,A) \to H^0 (X , A) \to H^0 (X \ohne K , A) \to H^1_K (X,A) \to H^1 (X,A) \; ,
\]
we obtain the exact sequence:
\[
0 \to H^0_c (X ,A) \to H^0 (X , A) \to C^{\infty} (E [X],A) \to H^1_c (X,A) \to H^1 (X,A) \; .
\]
On the other hand, let $j : X \hookrightarrow \oX$ be a compactification of $X$ and let $i : X_{\infty} = \oX \ohne X \hookrightarrow \oX$ denote the closed immersion of the complement. The distinguished triangle on $\oX$
\[
j_! A \longrightarrow R j_* A \longrightarrow i_* i^{-1} Rj_* A \longrightarrow \ldots
\]
gives the exact sequence
\[
 0 \to H^0_c (X ,A) \to H^0 (X,A) \to H^0 (X_{\infty} , i^{-1} Rj_* A) \to H^1_c (X ,A) \to H^1 (X,A) \; .
\]
There is a natural map from this sequence to the one above. The 5-lemma therefore gives us a canonical isomorphism
\begin{equation}
  \label{eq:31}
  C^{\infty} (E [X] , A) = H^0 (X_{\infty} , i^{-1} Rj_* A) \; .
\end{equation}
Now let us look at analogies in the \'etale topology of arithmetic varieties. Let $\eX$ be a proper $\Z$-scheme and define $\oeX = \eX \amalg \eX_{\infty}$ as in the beginning of \S\,2. Using formulas (2) and (9) we find
\begin{eqnarray*}
  H^0 (\eX_{\infty} , i^* Rj_* (\Z / l^n)) & = & H^0 (\eX^{\ann}_{\R} , \alpha^* (\Z / l^n)) \\
& = & H^0_{\et} (\eX \otimes \R , \Z / l^n)\\
& = & \map (\pi_0 (\eX \otimes \R) , \Z / l^n) \; .
\end{eqnarray*}
Comparing this with (21) suggests to view $\pi_0 (\eX \otimes \R)$ as a replacement for the space of ends of $\eX$. In particular the space of ends of $\spec \eo_k$ would be the finite set of archimedian valuations of $k$.
\section{The role of the infinite primes in the conjectural dynamics}
\label{sec:2}

In \cite{D2} we introduced a certain module $\Fh_{\ep} (M)$ attached to a motive $M$ over an algebraic number field $k$ and a place $\ep$ of $k$. It is important for expressing the local $L$-factor of $M$ at $\ep$ as a zeta-regularized determinant. For finite primes we interpreted these modules geometrically within a certain conjectural dynamical framework, \cite{D4} 3.22.

The present section is devoted to an analogous interpretation of $\Fh_{\ep} (M)$ for the infinite places $\ep$. This section is quite speculative. We assume that the reader is acquainted with the ideas of \cite{D5}.

In \cite{D5} \S\,5 it was argued that a motive $M$ over $k$ should give rise to a certain sheaf of $\Rh$-modules $\Fh (M)$ on $X_{\spec \eo_k}$ with an action
\[
\psi^t : (\phi^t)^{-1} \Fh (M) \to \Fh (M) \; .
\]
Here $\Rh$ is the sheaf of smooth real valued functions on $X_{\spec \eo_k}$ which are locally constant on the leaves of the one-codimensional foliation $\Bh$ of $X_{\spec \eo_k}$. There should even be a sheaf $\Fh (M)$ (or possibly a complex of sheaves) on $X_{\overline{\spec \eo_k}}$ prolonging the sheaf $\Fh (M)$ on $X_{\spec \eo_k}$. Because of the discussion in \cite{D3} \S\S\,5, 9 we think of the prolongation as some intermediate direct image
\begin{equation}
  \label{42}
  \Fh (M) = j_{!*} (\Fh (M) /_{X_{\spec \eo_k}}) \; .
\end{equation}
Here $j : X_{\spec \eo_k} \hookrightarrow X_{\overline{\spec \eo_k}}$ is the map induced by the inclusion $j : \spec \eo_k \hookrightarrow \overline{\spec \eo_k}$ and a theory of perverse sheaves of $\Rh$-modules with respect to the middle perversity is required. 

It appears that infinite primes $\ep$ of $k$ give rise to stationary points $x_{\ep}$ on $X_{\overline{\spec \eo_k}}$ such that $\kappa_{x_{\ep}} = \kappa_{\ep}$. 

Here $\kappa_{\ep} = -1$ if $\ep$ is complex and $\kappa_{\ep} = -2$ if $\ep$ is real. 

The real number $\kappa_{x_{\ep}}$ is defined by the formula $T_{x_{\ep}} \phi^t = e^{t\kappa_{x_{\ep}}}$ on the $1$-dimensional real vector space $T_{x_{\ep}} (X_{\overline{\spec \eo_k}}) / T_{x_{\ep}} \Bh$.

The ``stalk'' of the conjectural sheaf $\Fh (M)$ at an infinite place $\ep$ of $k$ is known by the considerations of \cite{D2} \S\,6, \cite{D6} \S\,3. It is given by
\[
\Fh_{\ep} (M) = \Gamma (\R , \xi^{\infty}_{\C} (M_{\ep} , \Vh^{\hullet} M_{\ep})) \quad \mbox{if $\ep$ is complex}
\]
and by
\[
\Fh_{\ep} (M) = \Gamma (\R^{\ge 0} , \xi^{\infty}_{\R} (M_{\ep} , \Vh^{\hullet} M_{\ep} , F_{\infty})) \quad \mbox{if $\ep$ is real} \; .
\]
Here $\xi^{\infty} = \xi^{\omega} \otimes_{\Ch^{\omega}} \Ch^{\infty}$ where $\xi^{\omega}$ is the Rees functor from filtered vector spaces (with involution) to locally free sheaves on the real analytic space $\R$ (resp. orbifold $\R^{\ge 0} = \R / \mu_2$) introduced in \cite{D6} \S\,2. On $\R$ resp. $\R^{\ge 0}$ a flow $\phi^t$ is given by $\phi^t_{\C} (r) = r e^{-t}$ resp. $\phi^t_{\R} (r) = re^{-2t}$. The sheaves $\xi^{\infty}$ carry an action $\psi^t : (\phi^t)^{-1} \xi^{\infty} \to \xi^{\infty}$. Together we get an action on its global sections and hence on $\Fh_{\ep} (M)$.

How to obtain $\Fh_{\ep} (M)$ defined as above from the sheaf $\Fh (M)$ on $X_{\overline{\spec \eo_k}}$? A natural idea is the following

 Let $\ep$ be a complex prime. Then there should exist a $\phi^t$-equivariant smooth embedding:
\begin{equation}
  \label{43}
  \iota_{\ep} : \R \hookrightarrow X_{\overline{\spec \eo_k}}
\end{equation}
such that $0 \in \R$ is mapped to $x_{\ep}$. Moreover we should have a $\psi^t$-equivariant isomorphism of sheaves of $\R$-vector spaces:
\begin{equation}
  \label{44}
  \iota^{-1}_{\ep} \Fh (M) \cong \xi^{\infty}_{\C} (M_{\ep} , \Vh^{\hullet} M_{\ep}) \; .
\end{equation}
Here $\psi^t$-equivariance means that
\[
\psi^t_{\xi^{\infty}_{\C}} \ent \iota^{-1}_{\ep} (\psi^t_{\Fh (M)})
\]
where the right hand map is defined as the composition:
\[
(\phi^t_{\C})^{-1} (\iota^{-1}_{\ep} \Fh (M)) = \iota^{-1}_{\ep} ((\phi^t)^{-1} \Fh (M)) \xrightarrow{\iota^{-1}_{\ep} (\psi^t)} \iota^{-1}_{\ep} \Fh (M)\; .
\]
Formula (\ref{44}) for $M = \Q (0)$ says that
\[
\iota^{-1}_{\ep} \Rh \cong \Ch^{\infty}_{\R} \; .
\]
Hence we are led to the conclusion that $\iota_{\ep}$ is transversal to the leaves. In particular the inverse image under $\iota_{\ep}$ of the leaf through the fixed point $x_{\ep}$ should be just the point $0 \in \R$. 

Another argument for transversality is the following: If $\iota_{\ep}$ is transversal to the leaves then it induces an isomorphism:
\[
T_0 \, \iota_{\ep} : T_0 \R \silo T_{x_{\ep}} X / T_{x_{\ep}} \Bh \; .
\]
Under this map the action of $T_{x_{\ep}} \phi^t$ on the right corresponds to the action of $T_0 \phi^t_{\C}$ on $T_0 \R$. Thus $e^{t\kappa_{x_{\ep}}} = e^{-t}$ i.e. $\kappa_{x_{\ep}} = -1 = \kappa_{\ep}$ as it must be. 

Let us investigate (23) further. As a dynamical system $\R$ decomposes into the three orbits $r < 0 , r = 0$ and $r > 0$. Set $x_{\pm} = \iota_{\ep} (\pm 1)$. Then by equivariance:
\[
\begin{array}{rcll}
\iota_{\ep} (r) & = & \phi^{-\log |r|} (x_-) & \mbox{for} \; r < 0 \\
\iota_{\ep} (0) & = & x_{\ep} \\
\iota_{\ep} (r) & = & \phi^{-\log r} (x_+) & \mbox{for} \; r > 0 \; .
\end{array}
\]
Moreover:
\[
\lim_{t\to + \infty} \phi^t (x_{\pm}) = x_{\ep}
\]
since e.g. 
\[
\lim_{t\to \infty} \phi^t (x_+) = \lim_{r \to 0 \atop r > 0} \phi^{-\log r} (x_+) = \lim_{r \to 0 \atop r > 0} \iota_{\ep} (r) = \iota_{\ep} (0) \; .
\]
Because of the argument at the beginning of \cite{D6} \S\,2 I do not think that $\iota_{\ep}$ extends to an equivariant embedding from $\C$ into $X$ whose pullback is $\xi^{\infty}$ over $\C$. This suggests that for every complex prime $\ep$ there should be exactly two different trajectories $\gamma^+$ and $\gamma^-$ transversal to the leaves flowing into $x_{\ep}$ for $t \to + \infty$. All embeddings $\iota_{\ep}$ are obtained by specifying a point on $\gamma_+ \cup \gamma_-$. Note here that $x_+$ and $x_-$ determine each other since the resulting embedding $\iota_{\ep}$ must be differentiable at $r = 0$.

It is easy to check that if (24) holds for $\iota_{\ep}$ corresponding to one choice of a point on $\gamma_+ \cup \gamma_-$ then it holds for the embedding $\iota_{\ep}$ corresponding to any other point on $\gamma_+ \cup \gamma_-$. We view $x_{\ep}$ or the unordered pair $\{ \gamma_+ , \gamma_- \}$ as the ``arithmetic point'' corresponding to $\ep$ and the different choices of $\iota_{\ep}$ i.e. of a point on $\gamma_+ \cup \gamma_-$ as the ``geometric points''. This complements the discussion on arithmetic and geometric points in $X$ corresponding to the finite primes in \cite{D4} 3.22.

Now let $\ep$ be a real prime. Then there should be a $\phi^t$-equivariant embedding:
\[
\iota_{\ep} : \R^{\ge 0} \hookrightarrow X \quad \mbox{with} \; \iota_{\ep} (0) = x_{\ep} \; .
\]
Here $\phi^t$ acts on $\R^{\ge 0}$ by $\phi^t_{\R} (r') = r' e^{-2t}$. The map $\iota_{\ep}$ should be smooth for the orbifold structure of $\R^{\ge 0} \cong \R / \mu_2$ given by the sheaf
\[
\Ch^{\infty}_{\R^{\ge 0}} := (sq_* \Ch^{\infty}_{\R})^{\mu_2} \; .
\]
Here $sq : \R \to \R^{\ge 0}$ is the squaring map $sq (r) = r^2$.
This means that pullback via $\iota_{\ep}$ maps $\iota^{-1}_{\ep} \Ch^{\infty}_X$ to $\Ch^{\infty}_{\R^{\ge 0}}$. In other words: Firstly, the restriction of $\iota_{\ep}$ to $(0, \infty) \subset \R^{\ge 0}$ should be smooth. Secondly, for any smooth function $F$ in a neighbourhood of $x_{\ep}$ the pullback function
\[
f (r') = (\iota^*_{\ep} F) (r') = F (\iota_{\ep} (r'))
\]
defined in some intervall $[0,\varepsilon)$ should have the property that $g (r) = f (r^2)$ is smooth in $(-\sqrt{\varepsilon} , \sqrt{\varepsilon})$. As above there should be a $\psi^t$-equivariant isomorphism of sheaves of $\R$-vector spaces:
\begin{equation}
  \label{45}
  \iota^{-1}_{\ep} \Fh (M) \cong \xi^{\infty}_{\R} (M_{\ep} , \Vh^{\hullet} M_{\ep} , F_{\infty}) \; .
\end{equation}
Again $\iota_{\ep}$ should be transversal to the leaves with the inverse image of the leaf through $x_{\ep}$ consisting only of $0 \in \R^{\ge 0}$. Transversality implies that the tangent map
\[
T_0 \, \iota_{\ep} : T_0 \, \R^{\ge 0} = (T_0 \R)^{\mu_2} \silo T_{x_{\ep}} X / T_{x_{\ep}} \Bh
\]
is an isomorphism. It transforms the action of $T_0 \phi^t_{\R}$ on $T_0 \R^{\ge 0}$ into the action of $T_{x_{\ep}} \phi^t$ on $T_{x_{\ep}} X / T_{x_{\ep}} \Bh$. Thus we find that $e^{-2t} = e^{\kappa_{x_{\ep}}t}$ i.e. $\kappa_{x_{\ep}} = -2 = \kappa_{\ep}$ as expected.

The system $(\R^{\ge 0}, \phi^t_{\R})$ has the two orbits $r'> 0$ and $r' = 0$. Set $x = \iota_{\ep} (1)$. Then
\[
\iota_{\ep} (r') = \phi^{-\halb \log r'} (x) \quad \mbox{for} \quad r' > 0 \quad \mbox{and} \quad \iota_{\ep} (0) = x_{\ep}
\]
and
\[
\lim_{t\to + \infty} \phi^t (x) = x_{\ep} \; .
\]
By a similar argument as before, for every real prime $\ep$ there should be exactly one trajectory $\gamma$ transversal to the leaves and flowing into $x_{\ep}$ for $t \to \infty$. The embeddings $\iota_{\ep}$ as above correspond to points $x$ on this trajectory $\gamma$. We view $x_{\ep}$ or the trajectory $\gamma$ as the ``arithmetic point'' corresponding to $\ep$ and the embeddings $\iota_{\ep}$ i.e. the points on $\gamma$ as the corresponding ``geometric points''.

The preceeding discussion on the geometric meaning of the stalks $\Fh_{\ep} (M)$ for $\ep \tei \infty$ complements the analogous but easier discussion for $\ep \nmid \infty$ in \cite{D4} 3.22. There $\Fh_{\ep} (M)$ is isomorphic to the global sections of $\Fh (M)$ restricted to the periodic orbit $\gamma_{\ep}$ corresponding to the prime ideal $\ep$. 
\bigskip

\noindent {\bf Remark.} Choose an embedding
\[
\iota_{\infty} : \R / \mu_2 \cong \R^{\ge 0} \hookrightarrow X_{\overline{\spec \Z}}
\]
as above corresponding to the infinite prime $p = \infty$ and set $x = \iota_{\infty} (1)$. Thus
\begin{equation}
  \label{46}
  \iota_{\infty} [r] = \phi^{- \log r} (x) \; .
\end{equation}
The inverse image of the leaf through $x_{\infty}$ should consist only of $0 \in \R / \mu_2$. The inverse image of any other leaf via $\iota_{\infty}$ is a countable subset of \\
$(\R /\mu_2) \ohne \{ 0 \} = \R^* / \mu_2 $. Because of the formula
\[
\iota_{\infty} [rs] = \phi^{-\log r} (\iota_{\infty} (s))
\]
the multiplicative structure of $\R^* / \mu_2$ is related to the flow. In our context the only natural countable subsets of the group $\R^* / \mu_2$ are the translates $r \Q^* / \mu_2 \subset \R^* /\mu_2$. My guess is that these are exactly the inverse images of those leaves which do not pass through $x_{\infty}$. The space of leaves of $X_{\overline{\spec \Z}}$ would thus be
\begin{equation}
  \label{47}
  \R^* / \Q^* \cup \{ 0 \} = \R / \Q^* \; ,
\end{equation}
with the flow acting by multiplication with $e^{-t}$. 

The same conclusion can be reached by the following completely different argument: Assuming that $X_{\overline{\spec\Z}}$ is an $F$-system it is isomorphic to $M \times_{\Lambda} \R$ where $M$ is any leaf and $\Lambda$ is the group of periods c.f. \cite{D4} 3.12. By the argument in remark b) before 4.3 of \cite{D4} we must have $\Lambda = \log \Q^*_+$. For the convenience of the reader let us recall the argument: For any $F$-system $X$ the group of periods $\Lambda$ contains the lengths of the closed orbits \cite{D4} 3.11. Hence in case of $X = X_{\spec \Z}$ we have $\log p \in \Lambda$ for all primes $p$ and therefore $\log \Q^*_+ \subset \Lambda$.

On the other hand the canonical line bundle $\uR (1)$ on an $F$-system $X$ is given by the representation
\[
\exp \verk l : \pi_1 (X)^{\abb} \longrightarrow \R^*_+
\]
where
\[
l : \pi_1 (X)^{\abb} \longrightarrow \Lambda \subset \R
\]
is the period homomorphism. Thus $\uR (1)$ has a $\Q$-structure if and only if $\exp \verk l$ takes values in $\Q^*_+$ i.e. if $\Lambda := \Imm l \subset \log \Q^*_+$. Now for $X = X_{\spec \Z}$ the line bundle $\uR (1)$ does have a $\Q$-structure namely $F_{\Q} (M)$ for the motive $M = \Q (1)$. Hence for $X_{\spec \Z}$ we find that $\Lambda = \log \Q^*_+$. 

Thus the space of leaves of $X_{\spec \Z}$ is $\R / \log \Q^*_+$ with the flow acting by translation. This is isomorphic via the map $y \mapsto \exp (-y)$ to $\R^*_+ / \Q^*_+$ with the flow acting by multiplication with $e^{-t}$. The leaf through $x_{\infty}$ is mapped into itself by the flow since the point $x_{\infty}$ is preserved. In the space of leaves it therefore becomes a fixed point under the flow. Thus the space of leaves of $X_{\overline{\spec \Z}}$ is isomorphic to
\[
\R^*_+ / \Q^*_+ \cup \{ 0 \} \cong \R / \Q^*
\]
with the flow acting by $e^{-t}$.

In Connes' approach \cite{C} the phase space is $\A / \Q^*$ with $\R$ acting by multiplication by $e^{-t}$ in the $\infty$-coordinate. If we view $\A_f$ as the typical leaf and $\A_f / \Q^*$ as the leaf at infinity then the space of leaves is $\R^* / \Q^* \cup \{ 0 \} = \R / \Q^*$ as in our picture. However the discussion of fundamental groups for instance in \cite{D4} \S\,4 shows that our generic leaves must have a much more complicated topological structure than $\A_f$.

There is the interesting possibility that for $t \to - \infty$ the trajectory through $x$ becomes dense in $X_{\overline{\spec \Z}}$ or at least in $X_{\spec \Z}$. If one could guess the uniform structure induced on $\R / \mu_2$ from $X_{\overline{\spec \Z}}$ via $\iota_{\infty}$ one could obtain $X_{\overline{\spec \Z}}$ or at least $X_{\spec \Z}$ as a uniform space by completion. This looks somewhat similar to recovering the finite id\`eles of $\Q$ by completing $\Q^*$ with respect to the uniform structure induced by the embedding $\Q^* \hookrightarrow \A^*_f$. 

\noindent
Mathematisches Institut\\
Westf. Wilhelms-Universit\"at\\
Einsteinstr. 62\\
48149 M\"unster\\
Germany\\
deninge@math.uni-muenster.de

\end{document}